\documentclass[11pt,epsf]{article}
 \usepackage{amssymb}
 \usepackage{amsthm}
 \usepackage{mathrsfs}
 \usepackage{pifont}
 \usepackage{amsmath}
 \usepackage{latexsym}
 \usepackage{graphicx}
 \usepackage{epsfig}
  \newcommand{\integ}[2]{\displaystyle \int_{#1}^{#2}}

  \newcommand{\no}{\noindent}

 \newcommand{\PP}{I\!\!P}

 \newcommand{\ind}{1\!\!1}

 \newcommand{\beqar}{\begin{eqnarray}}
 \newcommand{\eeqar}{\end{eqnarray}}

 \newcommand{\gb}{\beta}
 
 \def \cadlag {{c\`adl\`ag}~}
 \newcommand{\1}{1\!\!1}
 \newtheorem{definition}{Definition}[section]
 \newtheorem{thm}{Theorem}[section]
 \newtheorem{propo}{Proposition}[section]
 
 \newtheorem{remark}{Remark}[section]
 \newtheorem{corollary}{Corollary}[section]
 
  \def \R{\mathbb{R}}

 \def \P{\mathbb{P}}
 \def \B{\mathbb{B}}

 \def \E{\mathbb{E}}
 \def \V{\mathscr{V}}

 \def \T{\mathcal{T}}
 \def \cD{\mathcal{D}}
 \def \H{\mathcal{H}}
 \def \F{\mathcal{F}}
 \def \A{\mathcal{A}}

 \def \bf{\textbf}
 \def \it{\textit}
 \def \S{\mathcal{S}}
 \def \ds{\displaystyle}
 
\begin{document}
 \title{Stochastic Impulse Control of Non-Markovian Processes}
 \author{Boualem Djehiche \thanks{Department of Mathematics, The Royal Institute of Technology, S-100 44 Stockholm,
Sweden.\@ e-mail: boualem@math.kth.se}\,\,\,\,Said
Hamad\`ene\thanks{Universit\'e du Maine, D\'epartement de
Math\'ematiques, Equipe Statistique et Processus, Avenue Olivier
Messiaen, 72085 Le Mans, Cedex 9, France. e-mail:
ibtissam.hdhiri@univ-lemans.fr} \,\, and \,\,  Ibtissam
Hdhiri\thanks{Universit\'e du Maine, D\'epartement de
Math\'ematiques, Equipe Statistique et Processus, Avenue Olivier
Messiaen, 72085 Le Mans, Cedex 9, France. e-mail:
hamadene@univ-lemans.fr}}
\date{\today}
 \maketitle
 \begin{abstract}
 We consider a class of stochastic impulse control problems of
 general stochastic processes i.e. not necessarily Markovian.
Under fairly general conditions we establish
existence of an optimal impulse control. We also prove existence of combined optimal
 stochastic and impulse control of a fairly general class of diffusions with random coefficients. Unlike, in the
 Markovian framework, we cannot apply quasi-variational inequalities techniques. We rather derive the main results using
 techniques involving reflected BSDEs and the Snell envelope.
 \end{abstract}

 \medskip

 \bigskip AMS subject Classifications: 60G40; 60H10; 62L15; 93E20;
 49N25.
 \newline
 \medskip \newline
 \textbf{Keywords:} stochastic impulse control; Snell envelope; Stochastic
 control; backward stochastic differential equations; optimal stopping
 time.
 \date{}

 \bigskip

 \section{Introduction} Finding a stochastic impulse control policy amounts to determining the
 sequence of random dates at which the policy is exercised and the sequence of impulses describing the magnitude
 of the applied policies, which maximizes a given reward function. Given the general applicability of stochastic impulse
 control models in various fields such as finance, e.g. cash management (see Korn (1999) for an excellent survey and the
 textbook by Jeanblanc \it{et al.} (2005) and the references therein), and management of renewable
 resources (see e.g. Alvarez (2004), Alvarez and Koskel (2007) and the references therein), it is not
 surprising that the mathematical framework of such problems is well established (see Lepeltier-Marchal (1984),
 \O ksendal and Sulem (2006) and the references therein and the seminal textbook by Bensoussan and Lions (1984)
 on quasi-variational inequalities and impulse control). Indeed, in most cases, the impulse control problem is studied
 relying on quasi-variational inequalities, which is possible only through tacitly assuming that the underlying dynamics
 of the controlled system is Markovian and the instantaneous part of the reward function a
 deterministic function of the value of the process at a certain instant. These assumptions are
 obviously not realistic in most applications, such as in certain models in commodities trading. Even
 if the underlying process is Markov, the instantaneous part of the reward function may
 depend on the whole path of the process or is simply random.

In this study we consider a class of stochastic impulse control
problems where the underlying dynamics of the controlled system is
typically not Markov and where the instantaneous reward functional
is random, in which case, we cannot rely on the well established
quasi-variational inequalities technique to solve it. Instead, we
solve the problem using techniques involving reflected BSDEs and the
Snell envelope that seem suit well this general situation. The main
idea is to express the value-process of the control problem as a
Snell envelope and show that it solves a reflected BSDE, whose
existence and uniqueness are guaranteed provided some mild
integrability conditions of the involved coefficients.  This is done
through an appropriate approximation scheme of the system of
reflected BSDEs that is shown to converge to our value process. The
underlying approximating sequence is shown to be the value process
of an impulse control over strategies which have only a bounded
number of impulses, for which an optimal policy is also shown to
exist. Finally, passing to the limit, letting the number of impulses
become large, we prove existence of an optimal policy of our
stochastic impulse control problem.

The paper is organized as follows. In Section 2  we recall the main
tools on reflected BSDEs and Snell envelope we will use to establish
the main results. In Section 3, we formulate the considered
stochastic impulse control. In Section 4, we consider an appropriate
approximation scheme of the system of reflected BSDEs that is shown
to converge to our value process. In Section 5, we establish
existence of an optimal impulse control over strategies with a
bounded number of impulses, in Section 6, we prove existence of an
optimal impulse control over all admissible strategies. Moreover,
the corresponding value process is the limit of the sequence of
value processes associated with the optimal impulse control over
finite strategies, as their number becomes large. Finally, in
Section 7, we consider a mixed stochastic control and impulse
control problem of a fairly large class of diffusion processes that
are not necessarily Markovian. Using a Bene\v s-type selection
theorem, we derive an optimal policy using similar tools.

\section{Preliminaries and notation}
 \noindent Throughout this paper $(\Omega ,\mathcal{F},\PP)$ is a
 fixed probability space on which is defined a standard
 $d$-dimensional Brownian motion $\B=(\B_{t})_{0\le t\leq T}$ whose natural filtration is $%
 (\F_{t}^{0}:=\sigma \{\B_{s},s\leq t\})_{0\le t\leq T}$ ;
 $(\F_{t})_{0\le t\leq T}$ is the
 completed filtration of $(\F_{t}^{0})_{0\le t\leq T}$ with the $\PP$-null sets of $%
 \mathcal{F}$, hence $(\F_{t})_{0\le t\leq T}$ satisfies the usual
 conditions, $i.e.$, it is right continuous and complete. Let
 \begin{itemize}
 \item $\mathscr{P}$ be the $\sigma$-algebra on $[0,T]\times
 \Omega$ of $\F_t$-progressively measurable processes.
 \item  for any $p \leq 2$, ${\H}^{p,k}$  be the set of
 $\mathscr{P}$-measurable processes $v=(v_t)_{0\le t\leq T}$  with
 values in $\R^k$ such that $\E[\int_{0}^{T}|v_s|^pds]<\infty$.
\item ${\cal S}^2$ (resp. $\S_c^2$) be the set
 of $\mathscr{P}$-measurable and \cadlag (abbreviation of right continuous and left limited) (resp. continuous)
 processes $Y=(Y_t)_{0\le t\leq T}$ such that $\E[\sup_{0\le t\leq T}|Y_t|^2]<\infty$.
  \item ${\S}_i^2$ (resp. $\S_{c,i}^2$) the set of non-decreasing processes
    $k=(k_t)_{0\le t \leq T}$  of
  $\S^2$ (resp. $\S_c^2$) which satisfy $k_0=0$.
\item for $t\leq T$, ${\cal T}_t$   the set of $\F_t$-stopping times $\nu$ such
 that $\PP-a.s.,\,\,t \leq \nu \leq T$. Finally for any stopping time $\nu$, ${\cal F}_\nu$ is the
 $\sigma$-algebra on $\Omega$ which contains the sets $A$ of ${\cal F}$ such that $A \cap
\{\nu\leq t\}\in {\cal F}_t.$
 $\Box$
 \end{itemize}
 \medskip

 Consider now an ${\cal S}^2$-process  $X=(X_t)_{0\le t\leq T}$. The
 $Snell$ $envelope$ of $X$, which we denote by $N(X)=(N(X)_t)_{0\le t\leq
 T}$, is defined as
 $$
 \PP-a.s. \,\,\,\; N(X)_t=\mbox{ess sup}_{\nu \in {\cal
 T}_t}\E[X_\nu | \F_t], \,\,\,\, 0\le t\le T.
 $$
 It is the smallest \cadlag $(\F_t,\PP)$-supermartingale of
 class $[D]$ (see the appendix for the definition) which dominates $X$, $i.e.,$ $\PP-a.s.,\;N(X)_t\geq
 X_t$, for all $ 0\le t\leq T$.

\medskip For the sequel, we need the following result related to the
 continuity of the Snell envelope  with respect to increasing
 sequences whose proof can be found in Cvitanic and Karatzas (1996) or Hamad\`ene and Hdiri (2007).
 \begin{propo}\label{lem1} Let $(U_{n})_{n\geq 1}$ be a sequence of \cadlag and
 uniformly square integrable processes which converges increasingly
 and pointwisely to a \cadlag and uniformly square integrable
 process $U$, then $(N(U_{n}))_{n\geq 1}$ converges increasingly
 and pointwisely to $N(U)$.
 \end{propo}

\medskip\noindent In the Appendix at the end of the paper, we collect further results on the Snell envelope we
will refer to in the rest of the paper.

\medskip \noindent Let us underline that in the Markovian case, the problem
under consideration is solved using  PDEs techniques. However, in our
framework, we can no longer apply these techniques. Instead, we use backward
stochastic differential equations (BSDEs in short) which we will introduce with others properties.

\medskip
Let $X = (X_t)_{0\le t \leq T}$ be a barrier process of
 ${\cal S}^2$ and $f:[0,T]\times \Omega \times \R^{1+d}\mapsto \R$  a
 drift coefficient such that $(f(t,\omega,0,0))_{0\le t \leq T} \in
 {\H}^{2,1}$ and uniformly Lipschitz in $(y,z)$, i.e. there exists a constant
 $C>0$ such that
 $$
 |f(t,y,z) - f(t,y',z')| \leq C (|y-y'|+|z-z'|) \mbox{ for any }
 t,y,z,y' \mbox{ and }z'.
 $$
Then we have the following
 \begin{thm} $($Hamad\`ene $($2002$)$$).$ \label{thexist}
  There exists  a unique $\mathscr{P}$-measurable triple of
  processes $(Y,Z,K)=(Y_t,Z_t,K_t)_{0\le t\leq T}$ with values in $\R^{1+d+1}$
 solution of the reflected BSDE associated with $(f,X)$, i.e.,
 $$
 \left\{
 \begin{array}{ll}
     Y \in {\cal S}^2,\:Z\in {\H}^{2,d}\hbox{ and }  K \in \S_i^2,  \\
      Y_t=X_T+\integ{t}{T}f(s,Y_s,Z_s)ds +
 K_T-K_t-\integ{t}{T}Z_s d\B_s,\,\,\, 0\le t\leq T,\\
     Y_t\geq X_t,\;\,\,\,\text{for all}\,\,\, 0\le t\leq T, \\
     \int_0^T(Y_t-X_t)dK^c_t=0, \hbox{ and
     }\Delta_tY:=Y_t-Y_{t-}=-(X_{t-}-Y_t)^+\;\ind_{[Y_t-Y_{t-}<0]},
 \end{array}
 \right.
 $$
 where $K^c$ is the continuous part of $K$. Moreover, $Y$ admits
   the following representation.
 \begin{equation}\label{snell}
  \PP-a.s.,\;\; Y_t = \mbox{ess sup}_{\tau
 \in \T_t} \E[\int_t^\tau f(s,Y_s,Z_s)ds + X_\tau|\F_t],\,\,\, t\le T.
 \end{equation}
  \no In addition, if $X$ is left upper semi-continuous, $i.e.$,
 it has only positive jumps, then the process $Y$ is continuous.
\end{thm}

\medskip\noindent From (\ref{snell}) we note that  $(Y_t + \int_0^t
 f(s,Y_s,Z_s)ds)_{0\le  t \leq T}$ is the Snell envelope of the process
 $(\int_0^t f(s,Y_s,Z_s)ds + X_t)_{0\le t \leq T}.$

 \medskip\noindent In view of the results in El-Karoui \it{et al.} (1995), solutions of BSDEs with one reflecting barrier
 can be compared when we can compare the generators, the terminal
 values and the barriers. This remains true in this framework of
 discontinuous processes. Indeed, the following result holds.

\begin{propo} $($Hamad\`ene $($2002$)$$)$ \label{thcompa}
 Let $\widetilde{f}$ (resp. $\widetilde{X}$) be another
 map  from $[0,T]\times \Omega\times \R^{1+d}$ into $\R$ (resp.
 another process of ${\cal S}^2$) such that:
\begin{itemize}
\item[$(i)$] there exists a process
 $(\widetilde{Y},\widetilde{Z},\widetilde{K})=(\widetilde{Y}_t,\widetilde{Z}_t,\widetilde{K}_t)_{t\leq
   T}$ solution of the reflected BSDE associated with $(\widetilde{f},\widetilde{X})$
\item[$(ii)$] $\PP-a.s.$ $\forall t\leq T$, $f(t,
 \widetilde{Y}_t,\widetilde{Z}_t)\leq
 \widetilde{f}(t,\widetilde{Y}_t,\widetilde{Z}_t)$

\item[$(iii)$] $\PP-a.s.$,  for all $ t\leq T$, $X_t\leq
  \widetilde{X}_t$.
\end{itemize}
\no Then, we have $\PP-a.s.$,  for all $ t\leq T$, $Y_t\leq
 \widetilde{Y}_t$. $\Box$
\end{propo}

\no Now, let us consider a sequence $(y^n,z^n,k^n)_{n \geq 1}$ of
 processes defined as follows:
 \begin{equation*}
 \left\{
 \begin{array}{l}
 (y^n,z^n,k^n)\in {\S_c}^2\times {\H}^{2,d}\times {\S_{c,i}^2}, \\
 y^n_t=y_T^n+\integ{t}{T}f(s,y_s^n,z_s^n)ds +k^n_T - k^n_t
 -\integ{t}{T}
 z^n_s d\B_s, \,\, t\leq T,\\
 y^n_t \geq X_t,\;\, \text{for all}\,\, t\leq T, \mbox{ and }\integ{0}{T}(y^n_t
 -X_t)dk^n_t =0.
 \end{array}
 \right.
 \end{equation*}
 We now recall the following result by S. Peng (1999) which generalizes
 a well know property of supermartingales which tells that an
 increasing limit of \cadlag supermartingales is a also a \cadlag
 supermartingale.
 \medskip

 \begin{propo}$($Peng $($1999$, pp. 485)$$)$\label{propopeng}
 Assume the sequence $(y^n)_{n \geq 0}$ converges increasingly to a
 process $(y_t)_{0\le t \leq T}$ such that $\E[\sup_{0\le t \leq T} |y_t|^2]
 < \infty$  , then there exist two processes $(z,k) \in \H^{2,d}
 \times {\S_{i}}^2$ such that
 $$
 y_t=y_T+\int_t^T f(s,y_s,z_s) ds + k_T-k_t-\int_t^T z_s d\B_s.
 $$
 In addition, $z$ is the weak (resp. strong) limit of $z^n$ in
 $\H^{2,d}$ (resp. in $\H^{p,d}$, for $p < 2$) and for any stopping
 time $\tau$, the sequence $(k^n_\tau)_{n\geq
0}$ converges to $k_\tau$ in $L^p(dP)$. \end{propo}
 \noindent In this result, the assumption $\E[\sup_{0\le t \leq T} |y_t|^2] <
 \infty$ can be replaced by $\E[\sup_{n \geq 1} \sup_{t \leq T}
 |y_t^n|^2] < \infty$.

\section{ Formulation of the impulse control problem}
 Let  $L=(L_t)_{0\le t\leq T}$ be a stochastic process that describes the
 evolution of a system. We assume it $\mathscr{P}$-measurable, with values in $\R^l$
  and is such that $\E[\int_{0}^{T}|L_s|^2ds]<\infty$.
 An impulse control  is a sequence of  pairs $\delta = (\tau_n, \xi_n)_{n \geq 0}$ in which
 $(\tau_n)_{n \geq 0}$ is a sequence of ${\cal F}_t$-stopping times such that $0
 \leq \tau_0 \leq \tau_1\leq\ldots\leq T\;\:\PP$-a.s. and $(\xi_n)_{n
 \geq 0}$ a sequence of random variables with values in a \it{finite} subset $U$
 of $\R^l$  such that $\xi_n$ is ${\cal F}_{\tau_n}$-measurable. Considering
 the subset $U$ finite is in line with the fact that, in practice, the
 controller has only access to limited resources which allow him to exercise
impulses of finite size.

\medskip\noindent
The sequence $\delta = (\tau_n, \xi_n)_{n \geq 0}$ is said to
 be an admissible strategy of the control, and the set of admissible strategies will be denoted by ${\A}$.
 The controlled process $L^\delta=(L^\delta_t)_{0\leq t\leq T}$ is described as
 follows:
 \begin{equation}
 L_t^{\delta}=\left\{
 \begin{array}{ll} L_t,\,\,\,\text{if}\,\,\,0\leq t\leq \tau_0,
 \\ L_t+\xi_n,\,\,\,\text{if}\,\,\,\tau_n\leq t< \tau_{n+1},\,\,\,n\ge 0,
 \end{array}
 \right.
 \end{equation}
 or, in compact form,
 $$
 L^\delta_t=L_t+\sum_{n\geq 0}\xi_n\ind_{[\tau_n\leq
 t]},\,\,\, 0\le t\leq T.
 $$
The associated reward of controlling the system is
 $$
 J(\delta)=\E[\integ{0}{T}h(s,\omega,L^\delta_s)ds-\sum_{n\geq
 0}\psi(\xi_n) \ind_{[\tau_n<T]}],
 $$
 where $h$, represents the instantaneous reward and $\psi$ the  costs
 due to the impulses.

\medskip\noindent
This formulation of impulse control also falls within the class of singular
stochastic control problems, since the bounded variation part of the process,
 which controls the dynamic of the system, is allowed to be only purely discontinuous- See \O ksendal and Sulem (2006) for
 further details. Finally, note that if for example the process $L$
 satisfies
 $$
L_t=L_0+\int_0^tb(s,\omega)ds+\int_0^t\sigma(s,\omega)dB_s,\,\,\,\,
   t\leq T,
$$
where, $(b(s))_{0\le s\leq T}$ and $(\sigma(s))_{0\le s\leq T}$ are adapted stochastic processes, the existing theory on impulse control
 cannot be applied to the associated problem, since the processes $b$
 and $\sigma$ are random.

\medskip
We make the following assumptions  on $h$ and $\psi$.

\medskip\noindent \bf{Assumption (A)}
\begin{itemize}
\item [(A1)] $h:\,\,\, [0,T] \times \Omega \times \R^l\longrightarrow [0,+\infty)$  is uniformly bounded by a constant $\gamma$ in all its arguments i.e.\@ for any $ (t,\omega, x)\in [0,T] \times
 \Omega \times \R^l,\,\, 0\leq h(t,\omega,x)\leq \gamma$.
\item [(A2)] $\psi:\,\,\, U\longrightarrow [0,+\infty)$ is bounded from below, i.e. there exists a constant $c>0$ such that $\inf_{\beta\in U}\psi(\beta)\ge c$.
\end{itemize}

\medskip\noindent
Assumption (A2) is motivated by the following form of
proportional and fixed transaction costs (see Korn (1999) or Baccarin and
Sanfelici (2006) for further examples).
$$
\psi(\xi)=\phi(\xi)+c,
$$
where $\phi\geq 0$, $\phi(0)=0$ and $c$ is positive constant.
\medskip

\begin{definition}A strategy $\delta^*\in {\A}$ such that
\begin{equation}\label{opt.strategy}
J(\delta^*)=\sup_{\delta \in {\A}}J(\delta)
\end{equation}
is called optimal.
\end{definition}

\medskip\noindent The properties of $h$ and $\psi$ make the supremum of the reward function $J$ over the set $\A$
coincides with the one over the set of finite strategies, $\cD$ defined as
$$
\cD=\{ \delta=(\tau_n;\;\beta_n)_{n \geq 0}\in \A;\quad  \PP (\tau_n (\omega) <T,\,\,\, n \geq 0)= 0 \}.
$$
That is,
$$
 \sup_{\delta \in \A} J(\delta)=\sup_{\delta \in \cD} J(\delta).
$$
Indeed, consider a strategy  $\delta =(\tau_n;\;\beta_n)_{n \geq 0}$  of $\A$ which does not belong to $\cD$ and let $B = \{\omega \in \Omega;\,\, \tau_n (\omega) <T, \,\,\, n \geq 0 \}$. Since $\delta$ is not finite, $\PP(B)>0$. But, since $h$ is bounded, we have
\begin{eqnarray*}
 J(\delta) &=& \E[ \int_0^T h(s,L_s^\delta) ds-\sum_{n \geq 0}
\psi(\beta_n) \ind_{[\tau_n <T]}] \\
 &\leq& \gamma T-\E[(\sum_{n \geq 0} \psi(\beta_n)
 \ind_{[\tau_n <T]}) \ind_B-(\sum_{n \geq 0}\psi(\beta_n))
 \ind_{[\tau_n <T]}) \ind_{B^c}]\\ &=& -\infty,
 \end{eqnarray*}
 whence the desired result.

\section{An approximation scheme}

 For any stopping time $\nu $ and an ${\cal
 F}_{\nu}-$measurable random variable $\xi$, let $(Y^0_t(\nu,\xi),Z^0_t(\nu,\xi))_{0\le t\leq T}$ be the
 solution in ${\S_c}^2\times {\H}^{2,d}$ of the following standard
 BSDE :
 \begin{equation} \label{y0def}
 Y^0_t(\nu,\xi)=\integ{t}{T}h(s,L_s+\xi)\ind_{[s\geq \nu]}ds
 -\integ{t}{T} Z^0_s(\nu,\xi)d\B_s, \,\,\, 0\le t\leq T.
 \end{equation}
 The solution of this BSDE exists and is unique by the well known
 Pardoux-Peng's Theorem (see Pardoux and Peng (1990))
 since the terminal value is null and the
 function $h$ is bounded. Next, for any  $n\geq 1$, let
 $(Y^n_t(\nu,\xi),K^n_t(\nu,\xi),Z^n_t(\nu,\xi))_{0\le t\leq T}$ be the
 sequence of processes defined recursively as solutions of
 reflected BSDEs in the following way:
 \begin{equation} \label{yndef}
 \left\{
 \begin{array}{l}
 (Y^n(\nu,\xi),Z^n(\nu,\xi),K^n(\nu,\xi))\in {\S_c}^2\times {\H}^{2,d}\times
{\S_{c,i}^2}, \\
 Y^n_t(\nu,\xi)=\integ{t}{T}h(s,L_s+\xi)\ind_{[s\geq \nu]}ds
 +K^n_T(\nu,\xi)- K^n_t(\nu,\xi)-\integ{t}{T}
 Z^n_s(\nu,\xi)d\B_s, \,\,\,0\leq t\leq T,\\
 Y^n_t(\nu,\xi)\geq O^n_t(\nu,\xi):=\max_{\gb \in
U}\{-\psi(\gb)+Y^{n-1}_t(\nu,\xi+\gb)\}, \,\,\, 0\leq t\leq T, \\
\integ{0}{T}(Y^n_t(\nu,\xi)-O^n_t(\nu,\xi))dK^n_t(\nu,\xi)=0.
 \end{array}
 \right.
 \end{equation}
 \begin{propo}\label{propo1n}
 For any $n \geq 0$, $\nu \in \T_0$ and any $\F_\nu$-measurable r.v. $\xi$, the triple
 \\$(Y^n(\nu,\xi),K^n(\nu,\xi),Z^n(\nu,\xi))$ of (\ref{yndef}) is well posed.
 \noindent Moreover, it satisfies the following properties.

 $(i)$ $\,\,\,\P-a.s.\,\,\,\,\, 0\leq Y_t^n(\nu,\xi)\leq Y_t^{n+1}(\nu,\xi),\quad 0\leq t \leq T$.

$(ii)$ $\,\,\,\P-a.s.\,\,\,\,\, Y_t^n(\nu,\xi)\leq \gamma
(T-t),\quad 0\leq t \leq T$.
  \end{propo}
 \noindent \emph{Proof}: We prove the result by induction on $n$. We first begin to show the
 well-posedness of  $(Y^n(\nu,\xi),K^n(\nu,\xi),Z^n(\nu,\xi))$ for any $n\geq 0$. As pointed out previously for $n=0$,
 for any stopping time $\nu$ and any  $\F_\nu$-measurable $r.v.$ $\xi$,
 the pair $(Y^0(\nu,\xi),Z^0(\nu,\xi))$ exists and belongs to
${\cal S}_c^2\times {\cal H}^{2,d}$. Suppose now for some $n\geq 1$,
for any stopping time $\nu$ and any
$\F_\nu$-measurable $r.v.$ $\xi$, the triplet
$(Y^n(\nu,\xi),K^n(\nu,\xi),Z^n(\nu,\xi))$ exists and belongs to
${\cal S}_c^2\times {\cal S}_{c,i}^2\times {\cal H}^{2,d}$.
Hence, thanks to the finitness of $U$,
 $(O^{n+1}_t(\nu,\xi))_{0\leq t\leq T}$ is  a continuous process and
 satisfies $O^{n+1}_T(\nu,\xi)\leq 0$. In view of Theorem \ref{thexist}, the triplet
 $(Y^{n+1}(\nu,\xi),K^{n+1}(\nu,\xi),Z^{n+1}(\nu,\xi))$ exists and belongs to
${\cal S}_c^2\times {\cal S}_{c,i}^2\times {\cal H}^{2,d}$. Thus, for
any $n\geq 0$, any stopping time $\nu$ and any
$\F_\nu$-measurable $r.v.$ $\xi$, the triplet
$(Y^n(\nu,\xi),K^n(\nu,\xi),Z^n(\nu,\xi))$ exists and belongs to
${\cal S}_c^2\times {\cal S}_{c,i}^2\times {\cal H}^{2,d}$.

Let us now show $(i)$ and $(ii)$. Once more we will use an induction argument. First writing $Y^0_t(\nu,\xi)$ as a conditional expectation w.r.t.
${\cal F}_t$ and taking into account of $0\leq h\leq \gamma$ we
obtain that $0\leq Y^0_t(\nu,\xi) \leq \gamma (T-t)$, for any
stopping time $\nu$ and any $\F_\nu$-measurable $r.v.$ $\xi$. Next,
as $K^1(\nu,\xi)$ is an increasing process then using standard
comparison result of solutions of BSDEs (see e.g. El-Karoui \it{et al.}
(1995)), we obtain $Y^0(\nu,\xi)\leq
 Y^1(\nu,\xi)$. Therefore, Properties $(i)$ and $(ii)$ hold for
 $n=0$. Suppose now that for some $n$, for any stopping time $\nu$ and any
 $\F_\nu$-measurable $r.v.$ $\xi$, $(i)$ and $(ii)$ hold. Then,
$O^{n+1}(\nu,\xi)\leq O^{n+2}(\nu,\xi)$ and then the
characterization (\ref{snell}) implies that $Y^{n+1}(\nu,\xi)\leq
Y^{n+2}(\nu,\xi)$. On the other hand, since, for any $\zeta \in {\cal
F}_\nu$, $Y^n(\nu,\zeta)\leq
 \gamma (T-t)$, it holds that
 $O^{n+1}_t(\nu,\xi)=\max_{\beta \in U}(-\psi(\gb)+
 Y^n_t(\nu,\xi+\beta))\leq \max_{\beta \in U}(-\psi(\gb)+\gamma
 (T-t))\leq \gamma (T-t),\,\,\, 0\leq t\leq T$.

\medskip\noindent
Now, once more by
 (\ref{snell}), we have, for any $n \geq 1$,
 \begin{equation}\label{yn1}\begin{array}{l}
 Y_t^{n+1} (\nu,\xi) = \mbox{ess sup}_{\tau \in \T_t} \E[\int_t^\tau
 h(s,L_s+\xi) ds + O_\tau^{n+1}(\nu,\xi) \ind_{[\tau <T]}|\F_t],\, t
 \leq T.\end{array}
 \end{equation}
 Therefore,
 $$
 Y_t^{n+1} \leq \mbox{ess sup}_{\tau \in \T_t}
 \E[\gamma (\tau -t) + \gamma (T-\tau)|\F_t] = \gamma (T-t)$$ and
this completes the proof of the claim. $\Box$
\medskip

In the next proposition we identify the limit process
$Y_t(\nu,\xi):=lim_{n\rightarrow
 \infty}Y^n_t(\nu,\xi)$ (which exists according to the last proposition) as a
Snell envelope. Note that, as a limit of a non-decreasing sequence
of continuous processes, $Y(\nu,\xi)$ is upper semi-continuous.
Moreover, it holds that
 \begin{equation} \label{posity}
0\leq Y_t(\nu,\xi)\le \gamma (T-t), \,\,\,\mbox{for all}\,\, t\le T,\,\,\, \mbox{
  and }\,\,\, Y_T(\nu,\xi)=0.\end{equation}
 Finally, once more thanks to the finitness of $U$, the sequence of processes $(O^n
 (\nu,\xi))_{n \geq 0}$ converges to $O(\nu,\xi)$ as $n \rightarrow
 \infty$, where, $O_t(\nu,\xi) := \max_{\beta \in U} \left[-\psi(\beta) + Y_t(\nu,\xi+\beta)\right]$), $0\leq t\leq T$.
  \begin{propo}\label{propo1}
\begin{itemize}
  \item[$(i)$] Let $\nu$ and $\nu'$ be two
 stopping times such that $\nu\leq \nu^{\prime}$ and $\xi$ an $\F_{\nu}$-
 measurable random variable, then it holds that $\P-a.s.,\,\,\,
 Y_t(\nu,\xi)=Y_t(\nu',\xi)$ for all $ t\ge \nu^{\prime}$.

 \item[$(ii)$]
  For any stopping time $\nu$ and $\F_{\nu}-$measurable random
variable $\xi$, the process
  $Y(\nu,\xi)$ is \cadlag and satisfies:
 \begin{equation}\label{s1}
  Y_t(\nu,\xi)=\mbox{ess sup}_{\tau \in \T_t}
 \E[\integ{t}{\tau}h(s,L_s+\xi)\:\ind_{[s\geq \nu]}ds \,+
 \,\ind_{[\tau <T]} \:O_\tau(\nu,\xi)|{\cal F}_t],\; t \leq T.
 \end{equation}
\end{itemize}
\end{propo}
$Proof$: $(i)$ We proceed by induction on $n$. We note that the
solution of the BSDE
$$Y^0_t(\nu,\xi)=\int_{t}^{T}h(s,L_s+\xi)\ind_{[s\geq \nu]} ds -
 \int_t^T Z_s(\nu,\xi) d\B_s $$ is unique. It follows that, for any $\xi \in
 \F_{\nu}$, $Y_t^0(\nu,\xi)=Y_t^0(\nu^{\prime},\xi)$ for any $t\geq \nu^{\prime}$.
 Assume now that the property holds true for some fixed $n$. Then
 $O_t^{n+1}(\nu,\xi)=O_t^{n+1}(\nu^{\prime},\xi),\forall t\geq \nu'$. Once more the uniqueness of the solution of (\ref{yndef}) yields
 $Y_t^{n+1}(\nu,\xi)=Y_t^{n+1}(\nu^{\prime},\xi)$, $\forall t\geq \nu^{\prime}$.
 Hence the property holds true for any $n \geq 0$ and the desired
 result is obtained by taking the limit as $n \rightarrow \infty$.

$(ii)$ The sequence of processes $\left((Y^n_t(\nu,\xi) + \int_0^t
h(s,L_s +
 \xi)ds)_{0\leq t\leq T}\right)_{n \geq 0}$ is of \cadlag
 supermartingales which converges increasingly and pointwisely to
 the process\\ $\left(Y_t (\nu,\xi)+ \int_0^t h(s,L_s + \xi)ds\right)_{0\le t\leq
 T}$. Therefore, according to Dellacherie and Meyer (1980, p. 86) and taking into account (\ref{posity}), the limit is also a
 \cadlag supermartingale. It follows that the process $Y(\nu,\xi)$ is
 also \cadlag. Next, the processes $O^n(\nu,\xi)$, $n\geq 1$, are \cadlag and
 converge increasingly to $O(\nu,\xi)$. The rest of the proof is a direct consequence of Proposition \ref{lem1}. $\Box$

\begin{remark} Propositions \ref{propo1n} and \ref{propo1} are generalizations of Corollaries 7.6 and 7.7 in \O ksendal and Sulem (2006).
\end{remark}

\section{Optimal impulse control over bounded strategies}
In this section we establish existence of an optimal impulse control over
the set of strategies which have only a bounded number of impulses.
Indeed, for fixed $n\ge 0$, let $\A_n$ be the following set of bounded strategies:
$$
\A_n=\{(\tau_m,\xi_m)_{m\geq 1}\in {\cal D},\mbox{ such that
}\tau_{n}=T, \PP-a.s.\}.
$$
Then, the following result, which is a generalizations of Theorem 7.2 in \O ksendal and Sulem (2006), holds.
\begin{propo}\label{rkyn}
 For $n \geq 1$, we have
\begin{equation}\label{opti-finite}
Y_0^n (0,0)=sup_{\delta \in \A_n} J(\delta).
\end{equation}
In addition, there exists a strategy $\delta_n^* \in \A_n$ which is
optimal, $i.e.,$
\begin{equation}\label{J-opti-finite}
J(\delta_n^*)=\sup_{\delta \in \A_n}J(\delta).
\end{equation}
\end{propo}

\no \emph{Proof}. Let $ \delta_n^* = (\tau_k^n,\beta_k^n)_{k \geq
 0}$ be the strategy defined as follows.
 \begin{eqnarray*}
  \tau_0^n = \inf\{s \geq 0;\; O^n_s (0,0) =
  Y^n_s(0,0)\} \wedge T,
 \end{eqnarray*}
and
 \begin{eqnarray}\label{eqjus}
 O_{\tau_0^n}^n(0,0):= \max_{\beta \in U} (-\psi(\beta) + Y^{n-1}_{\tau_0^n}
 (0,\beta))= \max_{\beta \in U} (-\psi(\beta) + Y^{n-1}_{\tau_0^n}
 (\tau_0^n,\beta))\\=-\psi(\beta_0^n) + Y^{n-1}_{\tau_0^n}
 ({\tau_0^n},\beta_0^n),\nonumber
 \end{eqnarray}
 and, for any $k \in \{1,\ldots,n-1\}$,
 \begin{eqnarray*}\tau_k^n &=& \inf\{s \geq \tau_{k-1}^n;\;
 O^{n-k}_s(\tau^n_{k-1},\beta_0^n+\ldots+\beta^n_{k-1}) =
  Y^{n-k}_s(\tau^n_{k-1},\beta_0^n+\ldots+\beta^n_{k-1})\} \wedge T,
 \end{eqnarray*}
 \begin{equation*}
  \mbox{and }O_{\tau^n_{k}}^{n-k}(\tau^n_{k-1},\beta_0^n+\ldots+\beta^n_{k-1})=-\psi(\beta_k^n) + Y^{n-k-1}_{\tau_k^n}
 ({\tau_k^n},\beta_0^n+\ldots+\beta_{k-1}^{n}+\beta_k^n).
 \end{equation*}
Note that in (\ref{eqjus}) we have taken into account the fact that
$Y^{n-1}_{\tau_0^n}
 (0,\beta))= Y^{n-1}_{\tau_0^n}
 (\tau_0^n,\beta))$. This equality is valid since $\beta$ is deterministic and
 thanks to the uniqueness of  the solutions of
 BSDEs (\ref{yndef}) which define $Y^{n-1}
 (0,\beta)$ and $Y^{n-1}(\tau_0^n,\beta)$ for $t\ge {\tau_0^n}.$ Finally, $ \tau_n^n =
  T$ and $ \beta_n^n \in U$ arbitrary. The choice of $\beta_n$
 is not very significant since there are no impulses at $T$. We will show that $\delta_n^*$ is an  optimal strategy.

 For any $k \leq n$, the random variables
 $\gb_k^n$ are $\F_{\tau_k^n}-$  measurable. Thanks to (\ref{snell}) and (\ref{yndef})
  we obtain
  $$
 Y^n_0(0,0) = \sup_{\tau \in \T} \E[\int_0^\tau h(s,L_s) ds
 +\ind_{[\tau < T]}O^n_\tau(0,0)].
  $$
 Moreover, since the process $O^n(0,0)$ is continuous and $O^n_T (0,0) \leq 0$,
 then the stopping time $\tau_0^n$ is optimal  after $0$.
 Therefore,
 \begin{eqnarray}\label{y00n}
 Y^n_0(0,0)= \E[\int_0^{\tau^n_0} h(s,L_s) ds +\ind_{[\tau_0^n <
 T]}O^n_{\tau_0^n}(0,0)].
 \end{eqnarray}
 Now, since for any $n\ge 1$,
 \begin{equation}\begin{array}{ll}
 O^n_{\tau_0^n}(0,0)&=\max_{\beta \in U}\{-\psi(\beta) +
 Y_{\tau_0^n}^{n-1} (0,\beta)\}=\max_{\beta \in U}\{-\psi(\beta) +
 Y_{\tau_0^n}^{n-1} (\tau_0^n,\beta)\}\\{}&
 =-\psi(\beta_0^n) + Y_{\tau_0^n}^{n-1}
 (\tau_0^n,\beta_0^n).\end{array}
 \end{equation} The second equality is valid since for any $\beta \in U$ we have $Y_{\tau_0^n}^{n-1}
 (0,\beta)=
 Y_{\tau_0^n}^{n-1} (\tau_0^n,\beta)$.

\medskip\noindent Then, it holds that
$$
Y^n_0(0,0)=\E[\int_0^{\tau^n_0} h(s,L_s) ds +\ind_{[\tau_0^n <
 T]}(-\psi(\beta_0^n) + Y_{\tau_0^n}^{n-1} (\tau_0^n,\beta_0^n))].
$$
 But, once again using (\ref{snell}) and (\ref{yndef}), we
 have
 $$
 Y_{\tau_0^n}^{n-1} (\tau_0^n,\beta_0^n) = \mbox{ess sup}_{\tau \in
 \T_{\tau_0^n}} \E[\int_{\tau_0^n}^\tau h(s,L_s+\beta_0^n) ds
 +\ind_{[\tau < T]}O^{n-1}_\tau(\tau_0^n,\beta_0^n)|\F_{\tau_0^n}].
 $$
and $\tau_1^n$ is an
 optimal stopping time after $\tau_0^n$. It yields that
 \begin{eqnarray*}
 Y_{\tau_0^n}^{n-1} (\tau_0^n,\beta_0^n)&=&
 \E[\int_{\tau^n_0}^{\tau^n_1} h(s,L_s+\beta_0^n) ds
 +\ind_{[\tau_1^n <
 T]}O^{n-1}_{\tau_1^n}(\tau_0^n,\beta_0^n)|\F_{\tau_0^n}]\\
 &=& \E[\int_{\tau^n_0}^{\tau^n_1} h(s,L_s+\beta_0^n) ds
 +\ind_{[\tau_1^n < T]} (-\psi(\beta_1^n) + Y_{\tau_1^n}^{n-2}
 (\tau_1^n,\beta_0^n+\beta_1^n))|\F_{\tau_0^n}].
 \end{eqnarray*}
 By combining the last equality and (\ref{y00n}) we get
 \begin{eqnarray*}
 Y^n_0(0,0)&=&  \E[\int_0^{\tau^n_0} h(s,L_s) ds
 +\int_{\tau^n_0}^{\tau^n_1} h(s,L_s+\beta_0^n) ds +\ind_{[\tau_0^n
 <
 T]} (-\psi(\beta_0^n)) \\
 &+&\ind_{[\tau_1^n < T]} (-\psi(\beta_1^n)) + \ind_{[\tau_1^n <
 T]} Y_{\tau_1^n}^{n-2} (\tau_1^n,\beta_0^n+\beta_1^n)],
 \end{eqnarray*}
 since $[\tau_1^n <T] \subset [\tau_0^n<T]$ and $\ind_{[\tau_0^n <
 T]} \int_{\tau^n_0}^{\tau^n_1} h(s,L_s+\beta_0^n) ds =
 \int_{\tau^n_0}^{\tau^n_1} h(s,L_s+\beta_0^n) ds$.\\
 Repeating this argument as many times as necessary yields
 \begin{eqnarray*}
 Y^n_0(0,0)&=&  \E[\int_0^{\tau^n_0} h(s,L_s) ds + \sum_{1\leq k
 \leq n-1} \int_{\tau^n_{k-1}}^{\tau^n_k}
 h(s,L_s+\beta_0^n+\ldots+\beta_{k-1}^n) ds \\&+&\sum_{0\leq k \leq
 n-1}\{\ind_{[\tau_k^n <
 T]} (-\psi(\beta_k^n))\} +\ind_{[\tau_{n-1}^n < T]} Y_{\tau_{n-1}^n}^{0}
 (\tau_{n-1}^n,\beta_0^n+\ldots+\beta_{n-1}^n)].
 \end{eqnarray*}
 But, according to (\ref{y0def}) we have
 $$
 Y_{\tau_{n-1}^n}^{0} (\tau_{n-1}^n,\beta_0^n+\ldots+\beta_{n-1}^n) =
 \E [\int_{\tau_{n-1}^n}^T h(s,L_s+\beta_0^n+\ldots+\beta_{n-1}^n)
 ds|\F_{\tau_{n-1}^n}].
 $$
 Therefore,
 \begin{eqnarray*}
 Y^n_0(0,0)&=&  \E[\int_0^{\tau^n_0} h(s,L_s) ds + \sum_{1\leq k
 \leq n} \int_{\tau^n_{k-1}}^{\tau^n_k}
 h(s,L_s+\beta_0^n+\ldots+\beta_{k-1}^n) ds +\sum_{0\leq k \leq
 n}\{\ind_{[\tau_k^n < T]} (-\psi(\beta_k^n))\}]\\
 &=&  \E[\int_0^{\tau^n_0} h(s,L_s) ds + \sum_{k \geq 1}
 \int_{\tau^n_{k-1}}^{\tau^n_k}
 h(s,L_s+\beta_0^n+\ldots+\beta_{k-1}^n) ds +\sum_{k \geq
 0}\{\ind_{[\tau_k^n < T]}
 (-\psi(\beta_k^n))\}]\\
 &=& J(\delta_n^*).
 \end{eqnarray*}

 \noindent It remains to show that $J(\delta_n^*) \geq J({\delta'}^n)$ for any
 strategy ${\delta'}^n$  of $\A_n$.
 \medskip\\
 Indeed, let ${\delta'}^n = ({\tau'}_0^n,\ldots,{\tau'}_{n-1}^n,T,T,\ldots;\;
 {\beta'}_0^n,\ldots,{\beta'}_{n-1}^n,{\beta'}_{n}^n,{\beta'}_{n}^n,\ldots)$
 be a strategy of $\A_n$. Since $\tau_0^n$ is optimal after $0$, we have
 \begin{eqnarray*}
 Y^n_0(0,0) &\geq& \E[\int_0^{{\tau'}_0^n} h(s,L_s) ds
 +\ind_{[{\tau'}_0^n < T]} O^n_{{\tau'}_0^n}(0,0)].
 \end{eqnarray*}
But,
$$
O^n_{{\tau'}_0^n}(0,0)=\max_{\beta\in U}\{-\psi(\beta)+Y^{n-1}_{{\tau'}_0^n}(0,\beta)\}
 = \max_{\beta\in U}\{-\psi(\beta)+Y^{n-1}_{{\tau'}_0^n}({\tau'}_0^n,\beta)\}
 \ge -\psi({\beta'}^n_0)+Y^{n-1}_{{\tau'}_0^n}({\tau'}_0^n,{\beta'}^n_0).
$$
Therefore, we have
 \begin{eqnarray*}
 Y^n_0(0,0) &\geq& \E[\int_0^{{\tau'}_0^n} h(s,L_s) ds
 +\ind_{[{\tau'}_0^n < T]}
 (-\psi({\beta'}_0^n) + Y_{{\tau'}_0^n}^{n-1} ({\tau'}_0^n,{\beta'}_0^n))]\\
 &\geq& \E[\int_0^{{\tau'}_0^n} h(s,L_s) ds
 +\int_{{\tau'}_0^n}^{{\tau'}_1^n} h(s,L_s+{\beta'}_0^n) ds
 +\ind_{[{\tau'}_0^n < T]} (-\psi({\beta'}_0^n)) \\
 &+& \ind_{[{\tau'}_1^n < T]} Y_{{\tau'}_1^n}^{n-2}
 ({{\tau'}_1^n},{\beta'}_0^n+{\beta'}_1^n)].
 \end{eqnarray*}
 Finally, iterating as many times as necessary we
 obtain
 \begin{eqnarray*}
 Y^n_0(0,0) &\geq& \E[\int_0^{{\tau'}^n_0} h(s,L_s) ds +
 \sum_{1\leq k \leq n} \int_{{\tau'}^n_{k-1}}^{{\tau'}^n_k}
 h(s,L_s+{\beta'}_0^n+\ldots+{\beta'}_{k-1}^n) ds \\&{}&+\sum_{0\leq k \leq
 n}\{\ind_{[{\tau'}_k^n < T]} (-\psi({\beta'}_k^n))\}]=J({\delta'}^n).
 \end{eqnarray*}
 Hence,  $J(\delta_n^*) \geq J({\delta'}^n)$, for any ${\delta'}^n
 \in \A_n$ The proof is now complete. $\Box$
\section{An optimal impulse control result.}
We now give the main result of this paper.
 \begin{thm}\label{veri-impulse} Under Assumption (A),
 the strategy $\delta^*=(\tau_n^*,\beta_{n}^*)_{n \geq 0}$ defined by
$$
\tau_0^* = inf\{s \geq 0;\; O_s (0,0)= Y_s(0,0)\}\wedge T,
$$
$$
\max_{\beta \in U} (-\psi(\beta) + Y_{\tau_{0}^*} (0,\beta)) =
-\psi(\beta_0^*) + Y_{\tau_{0}^*} (\tau_{0}^*,\beta_0^*),
$$
for $n \geq 1$,
$$
\tau_n^* = inf\{s \geq \tau_{n-1}^*;\; Y_s (\tau_{n-1}^*,\beta_0^*
+\ldots+\beta_{n-1}^*) =O_s (\tau_{n-1}^*,\beta_0^*
+\ldots+\beta_{n-1}^*)\}\wedge T,
$$
and
$$
 \max_{\beta \in U} (-c-\psi(\beta) + Y_{\tau_{n}^*}
(\tau_{n-1}^*,\beta_0^* +\ldots +\beta_{n-1}^*+\beta)) =
-c-\psi(\beta_0^*) + Y_{\tau_{n}^*} (\tau_{n}^*,\beta_0^*
+...+\beta_{n-1}^*+\beta_n^*).
$$
is optimal  for the impulse control problem.

\medskip Furthermore, we have
$$
Y_0(0,0)=J(\delta^*).
$$
 \end{thm}
\noindent\emph{Proof.} The proof is performed in three steps.

\noindent \bf{Step 1.} \it{Continuity of the value process
$(Y_t(\nu,\xi))_{0\le t\leq T}$}. We note that, by (\ref{s1}), we have, for any $0\leq t\leq
T$,
$$
Y_t(\nu,\xi)+\int_0^t h(s,L_s+\xi)\ind_{[s\geq \nu]}ds =
\mbox{ess sup}_{\tau \in \T_t}
\E[\integ{0}{\tau}h(s,L_s+\xi)\:\ind_{[s\geq \nu]}ds \,+
\,\ind_{[\tau <T]} \:O_\tau(\nu,\xi)|{\cal F}_t],
$$
meaning that the process $\left(Y_t(\nu,\xi)+\int_0^t h(s,L_s+\xi)\ind_{[s\geq
\nu]}ds\right)_{0\leq t \leq T}$ is the Snell envelope of \\
$\left(\int_{0}^{t}h(s,L_s+\xi)\:\ind_{[s\geq \nu]}ds \,+ \,\ind_{[t
<T]} \:O_t(\nu,\xi)\right)_{0\leq t \leq T}.$ Therefore, using
Proposition \ref{prapp2}, in the appendix below, there exist a
continuous martingale $M(\nu,\xi)$ and two increasing processes
$A(\nu,\xi)$ and $B(\nu,\xi)$ belonging to $\S_i^2$ such that
$B_0(\nu,\xi)=0$ and, for $0\leq t \leq T$,
$$
 \integ{0}{t}h(s,L_s+\xi)\:\ind_{[s\geq \nu]}ds
+Y_t(\nu,\xi) = M_t(\nu,\xi) - A_{t}(\nu,\xi) - B_t(\nu,\xi).
$$
In addition, the process $A(\nu,\xi)$ is optional and continuous,
and $B(\nu,\xi)$ is predictable and purely discontinuous. The
continuity of the value process $Y_t(\nu,\xi)$, will follow once
we show that, for any stopping time $\nu$ and $\F_\nu$-measurable
random variable $\xi$,  $B(\nu, \xi) \equiv 0$. Indeed, let us
assume that $B(\nu, \xi)$ is different to zero. Since the process is
non-decreasing and purely discontinuous, there exists $\tau \in
\T_\nu$ such that $B_\tau (\nu, \xi)-B_{\tau-} (\nu, \xi) >0$.
Thanks to (\ref{eqann}), in the appendix, we have $Y_{\tau-} (\nu,
\xi) = O_{\tau-} (\nu, \xi)$. Hence,
\begin{equation*}
Y_{\tau-} (\nu, \xi) = \max_{\beta \in U} (-\psi(\beta)+Y_{\tau-}
(\nu, \xi+\beta))>Y_{\tau} (\nu, \xi)\geq O_{\tau} (\nu,
\xi)=\max_{\beta \in U}(-\psi(\beta)+Y_{\tau} (\tau, \xi+\beta)).
\end{equation*}
Therefore, since $U$ is finite, there exists $\beta_1 \in U$ such that
the set
$$
\Lambda_1=\{Y_{\tau-} (\nu, \xi) = -\psi(\beta_1)+Y_{\tau-} (\nu,
\xi+\beta_1) \mbox{ and } \Delta Y_{\tau-} (\nu, \xi+\beta_1)<0\}$$
satisfies $P(\Lambda_1)>0$. But, the same holds for $\Delta
Y_{\tau-} (\nu, \xi+\beta_1)$. Therefore, there exists $\beta_2\in U$
such that the set
$$
\Lambda_2=\{Y_{\tau-} (\nu, \xi+\beta_1) = -\psi(\beta_2)+Y_{\tau-}
(\nu, \xi+\beta_1+\beta_2) \mbox{ and } \Delta Y_{\tau-} (\nu,
\xi+\beta_1+\beta_2)<0\}$$ satisfies $P[\Lambda_1 \cap
\Lambda_2]>0.$ It follows that, on the set $\Lambda_1 \cap \Lambda_2$,
we have
$$
Y_{\tau-} (\nu, \xi) = -\psi(\beta_1)-\psi(\beta_2)+Y_{\tau-} (\nu,
\xi+\beta_1+\beta_2) .
$$
Making this reasoning as many times as
necessary we obtain the existence of $\beta_1,\ldots,\beta_n$ elements
of $U$ and a subset $\Lambda_n$ of positive probability such that, on
$\Lambda_n$, we have
$$
Y_{\tau-} (\nu, \xi) = -\sum_{i=1}^n\psi(\beta_i)+Y_{\tau-} (\nu,
\xi+\beta_1+\ldots+\beta_n) \leq -nc+\gamma T.
$$
But, this is impossible
for $n$ large enough since the process $Y(\tau,\xi)$ is
non-negative. Therefore, the purely discontinuous process
$B(\nu,\xi)$ has no jumps and then it is null. Thus, the process
$Y(\nu,\xi)$ is continuous.
\medskip

\noindent \bf{Step 2.} The strategy
$\delta^*=(\tau_n^*,\beta^*_n)_{n\ge 0} \in \mathcal{D}$ and is such
that $Y_0(0,0) = J(\delta^*).$
\medskip

\no Using Proposition \ref{propo1}, we get
\begin{equation}\label{eqto} Y_0(0,0)= \mbox{ess sup}_{\tau \in \T_0}
\E[\integ{0}{\tau}h(s,L_s) ds \,+ \,\ind_{[\tau <T]} \:O_\tau(0,0)].
\end{equation}
Now, since $Y(\nu,\xi)$ is continuous for any $\nu \in \T$ and
any $\F_\nu$-measurable random variable $\xi$ and $O_T(0,0)\leq 0$,
then the stopping time $\tau_0^*$ is optimal for the problem
(\ref{eqto}). This yields
\begin{eqnarray}\label{y00}
Y_0(0,0)&=& \E[\integ{0}{\tau_0^*}h(s,L_s) ds \,+
\,\ind_{[{\tau_0^*} <T]} \:O_{\tau_0^*}(0,0)]. \nonumber
\end{eqnarray}
But,
$$
O_{\tau_0^*}(0,0)=\max_{\beta \in
U}\{-\psi(\beta)+Y_{\tau_0^*}(0,\beta)=\max_{\beta \in
U}\{-\psi(\beta)+Y_{\tau_0^*}(\tau_0^*,\beta)\}=-\psi(\beta^*)+Y_{\tau_0^*}(\tau_0^*,\beta^*)$$
where $\beta^*\in {\cal F}_{\tau_0^*}$. Note that the second
equality is valid thanks to Proposition \ref{propo1}-$(i)$.
Therefore,
$$
Y_0(0,0)= \E[\integ{0}{\tau_0^*}h(s,L_s) ds \,+ \,\ind_{[{\tau_0^*}
<T]} ](-\psi(\beta^*)+Y_{\tau_0^*}(\tau_0^*,\beta^*))].$$ Next,
\begin{eqnarray*}
Y_{\tau_0^*}(\tau_0^*,\beta_0^*)&=& \E[\integ{{\tau_0^*}}{\tau_1
^*}h(s,L_s+\beta_0^*) ds \,+ \,\ind_{[{\tau_1^*}
<T]} \:O_{\tau_1^*}(\tau_0^*,\beta_0^*)|\F_{\tau_0^*}]\\
&=& \E[\integ{\tau_0^*}{\tau_1^*}h(s,L_s+\beta_0^*) ds \,+
\,\ind_{[{\tau_1^*} <T]}
(-\psi(\beta_1^*)+\:Y_{\tau_1^*}(\tau_1^*,\beta_0^*+\beta_1^*))|\F_{\tau_0^*}].
\end{eqnarray*}
Replacing $Y_{\tau_0^*}(\tau_0^*,\beta_0^*)$ by its expression in
(\ref{y00}), we obtain
\begin{eqnarray*}
Y_0(0,0)&=&\E[\integ{0}{\tau_0^*}h(s,L_s) ds \,+
\,\integ{{\tau_0^*}}{\tau_1 ^*}h(s,L_s+\beta_0^*) ds
\,+\,(-\psi(\beta_0^*))\ind_{[{\tau_0^*}
<T]}\\&+&(-\psi(\beta_1^*))\,\ind_{[{\tau_1^*}
<T]}+\:Y_{\tau_1^*}(\tau_1^*,\beta_0^*+\beta_1^*)\,\ind_{[{\tau_1^*}
<T]}]
\end{eqnarray*} since $[\tau_1^*<T]\subset [\tau_0^*<T]$ and $[\tau_0^*<T]\in {\cal F}_{\tau_0^*}$. Proceeding in the same way as
many times as necessary we get
\begin{eqnarray}\label{y001}
Y_0(0,0)&=&\E[\integ{0}{\tau_0^*}h(s,L_s) ds +\ldots+
\,\integ{{\tau_{n-1}^*}}{\tau_n
^*}h(s,L_s+\beta_0^*+\ldots+\beta_{n-1}^*) ds
\,+\,(-\psi(\beta_0^*))\ind_{[{\tau_0^*}
<T]}+\ldots\nonumber\\&+&(-\psi(\beta_n^*))\,\ind_{[{\tau_n^*}
<T]}+\:Y_{\tau_n^*}(\tau_n^*,\beta_0^*+\ldots+\beta_{n-1}^*+\beta_n^*)\,\ind_{[{\tau_n^*}
<T]}].
\end{eqnarray}
Let us now show  that $\delta^* \in \mathcal{D}$. Assume that
$\P\{\tau_n^* < T;\;\; n \geq 0\}>0$. Then we have
\begin{eqnarray*}
Y_0 (0,0)&\leq& \E[\integ{0}{\tau_0^*}|h(s,L_s)| ds +\ldots\,+
\,\integ{{\tau_{n-1}^*}}{\tau_n
^*}|h(s,L_s+\beta_0^*+\ldots+\beta_{n-1}^*) |ds +\:\sup_{s \leq
T}|Y_s(\tau_n^*,\beta_0^*+\ldots+\beta_n^*)|\\
&+& \ind_{\{\tau_n^* < T;\;\; n \geq 0\}}\sum_{0\leq k \leq
n}(-\psi(\beta_k^*))\ind_{[{\tau_k^*} <T]}+\, \ind_{{\{\tau_n^* <
T;\;\; n \geq 0\}}^c}\sum_{0\leq k \leq
n}(-\psi(\beta_k^*))\ind_{[{\tau_k^*} <T]}]\\
&\leq& \gamma T+\:\E[\sup_{s \leq
T}|Y_s(\tau_n^*,\beta_0^*+\ldots+\beta_n^*)|]- nc\;\P{\{\tau_n^*
< T;\;\; n \geq 0\}}.
\end{eqnarray*}
The last quantity tends to $- \infty$ as $n \rightarrow \infty$,
then $Y_0(0,0) = - \infty$ which contradicts the fact that $Y(0,0)
\in \S^2$. Therefore, $\P\{\tau_n^* < T;\;\; n \geq 0\}=0$ i.e.
$\delta^* \in \mathcal{D}$. Finally, by taking limit as $n
\rightarrow \infty$ in (\ref{y001}) we obtain $Y_0(0,0) =
J(\delta^*)$.

%

\medskip\noindent \bf{Step 3.} \it{$J(\delta^*) \geq J(\delta)$ for any
  strategy $\delta \in \A$}. Let $\delta = (\tau_n,\beta_n)_{n \geq 0}$ be a finite strategy. Since $\tau_0^*$
 is optimal after $0$, we have
\begin{eqnarray*}
{Y}_0 (0,0)&\geq& \E[\int_0^{\tau_0} h(s,L_s) ds + \ind_{[\tau_0 <
T]} {O}_{\tau_0} (0,0)]\\&\geq& \E[\int_0^{\tau_0} h(s,L_s) ds +
\ind_{[\tau_0 < T]} \{-\psi(\beta_0)+ Y_{\tau_0}
(\tau_0,\beta_0)\}].
\end{eqnarray*}
But,
$$
{O}_{\tau_0} (0,0)=\max_{\beta \in
U}\{-\psi(\beta)+Y_{\tau_0}(0,\beta)\}=\max_{\beta \in
U}\{-\psi(\beta)+Y_{\tau_0}(\tau_0,\beta)\}\geq -\psi(\beta_0)+
Y_{\tau_0} (\tau_0,\beta_0).$$ It follows that
\begin{eqnarray*}
{Y}_0 (0,0) \ge \E[\int_0^{\tau_0} h(s,L_s) ds + \ind_{[\tau_0 < T]}
\{-\psi(\beta_0)+ Y_{\tau_0} (\tau_0,\beta_0)\}].
\end{eqnarray*}Next,
\begin{eqnarray*}
Y_{\tau_0} (\tau_0,\beta_0) &=& \mbox{ess sup}_{\tau \in
\T_{\tau_0}} \E[\int_{\tau_0}^{\tau} h(s,L_s+\beta_0)
ds+\ind_{[\tau < T]}
O_{\tau} (\tau_0,\beta_0)|\F_{\tau_0}]\\
&\geq& \E[\int_{\tau_0}^{\tau_1} h(s,L_s+\beta_0) ds+\ind_{[\tau_1 <
T]} \{-\psi(\beta_1)+ Y_{\tau_1}
(\tau_1,\beta_0+\beta_1)\}|\F_{\tau_0}].
\end{eqnarray*}
Therefore,
\begin{eqnarray*}
Y_0 (0,0) &\geq& \E[\int_0^{\tau_0} h(s,L_s)
ds+\int_{\tau_0}^{\tau_1} h(s,L_s+\beta_0) ds+
(-\psi(\beta_0))\ind_{[\tau_0 < T]}\\&+&
(-\psi(\beta_1))\ind_{[\tau_1 < T]}) +\ind_{[\tau_1 < T]}
Y_{\tau_1} (\tau_1,\beta_0+\beta_1)].
\end{eqnarray*}
Now, by following this reasoning as many times as necessary we
obtain,
\begin{eqnarray*}
Y_0 (0,0) &\geq& \E [\int_0^{\tau_0} h(s,L_s) ds+\sum_{1 \leq k
\leq n} \int_{\tau_{k-1}}^{\tau_k}
h(s,L_s+\beta_0+\ldots+\beta_{k-1}) ds
\\&+& \sum_{0 \leq k \leq n} (-\psi(\beta_k)) \ind_{[\tau_k <T]} +
 Y_{\tau_n} ({\tau_n},\beta_0+\ldots+\beta_n)]
\end{eqnarray*}
and since the strategy $\delta$ is finite, by taking the limit as $n
\rightarrow \infty$, we obtain $Y_0(0,0) \geq J(\delta)$ since
$|Y_{\tau_n} ({\tau_n},\beta_0+\ldots+\beta_n)|\leq \gamma
\1_{[\tau_n<T]}.$ As $\delta \in \A$ is arbitrary, then $Y_0(0,0) =
J(\delta^*)=\sup_{\delta \in \mathcal{D}} J(\delta)=\sup_{\delta \in
\A} J(\delta)$. $\quad\Box$

\begin{corollary} Under Assumptions (A) and (B) it holds that
\begin{equation}\label{limit}
\sup_{\delta\in
\mathcal{A}}J(\delta)=Y_0(0,0)=\lim_{n\to\infty}Y_0^n(0,0)=\lim_{n\to\infty}\sup_{\delta\in
\mathcal{A}_n}J(\delta).
\end{equation}
\end{corollary}

\section{Combined stochastic and impulse controls}
In this section we study a mixed stochastic and impulse control
problem, where, we allow the process $L$, that describes the
evolution of the system and subject to impulses, to also depend  on a
control $u$ from some appropriate set $\V$. Therefore, the dynamics
of the system is subject to a combination of control and impulses.
To begin with, we describe this dynamics.

\medskip\no
Let $\cal C$ be the set of continuous functions ${w}$ from $[0,T]$
into $\R^d$ endowed with the uniform norm. For $t\leq T$, let ${\cal
G}_t$ be the $\sigma$-field of $\cal C$ generated by
$\{\pi_s:w\mapsto w_s, \,\,s\leq t\}$. By $\mathbb{G}$ we denote the
$\sigma$-field on $[0,T]\times {\cal C}$ consisting of all the
subsets $G$, which have the property that the section of $G$ at time
$t$ is in ${\cal G}_t$ and the section of $G$ at $w$ is Lebesgue
measurable (see Elliott (1976) for more details on this subject).
Finally if $w\in \cal C$ and $a$ is a deterministic function then
$w+a$ is the function which with $t\in [0,T]$ associates
$(w+a)_t=w_t+a$.

\medskip\no
Let us now consider a function from $[0,T]\times {\cal
C}\rightarrow \R^d$ which satisfies the following
\medskip

\noindent \bf{Assumption (H)}.

\noindent (H1)  $\sigma$ is $\mathbb{G}$-measurable and there exists
a constant $k$ such that
\begin{itemize}
 \item [$(i)$] for every $t \in [0,T]$ and every $w$ and $ w'$ in $\cal{C}$,
$|\sigma(t,{w})-\sigma(t,{w'})| \leq k\|{w}-{w'}\|_t$ where
 $\|{w}\|_t = \sup_{s \leq t} |{w}_s|,\; t \leq T$;

\item [$(ii)$] for every $t \in [0,T]$, $|\sigma (t,0)| \leq k$, $\sigma$ is
 invertible and its inverse $\sigma^{-1}$ is bounded.
\end{itemize}

\medskip\no
Let $\cal V$ be a compact metric space and $\V$ the
 set of $\mathscr{P}-$measurable processes $v = (v_t)_{t \leq T}$
 with values in $\cal V$. Hereafter, $\V$ is called the set of admissible
 controls.

\medskip\no We consider now the process $(L_t)_{0\le t \leq T}$ which is the unique solution for the
 following stochastic differential equation:
 $$\left\{
 \begin{array}{ll}
              dL_t = \sigma(t,L_{\cdot})\, d\B_t,\;\,\, 0<t\leq T,\\
              L_0 = x,\;\,\,\, x\in \R^d,
 \end{array}
 \right.$$
whose existence is guaranteed by Assumption (H1). The process $L$ stands for the state of the
system when non-controlled.

 \medskip\no Let $f$ and (resp. $h$) be a  measurable and
 uniformly bounded function from $[0,T] \times {\cal C}  \times \cal
 V$ into $\R^d$ (resp. $\R^+$) such that
\begin{itemize}
\item [(H2)] $f$ and $h$ are $\mathbb{G}\otimes {\cal B}(\cal V)$-measurable

 \item [(H3)] for every $t\in [0,T]$, $w\in {\cal C}$, the function
 which with $u \in \cal V$ associates $f(t,w,u)$ (resp. $h(t,w,u)$) is
 continuous.
\end{itemize}
Now, given a control $u \in \V$, let  $\PP^u$ be the probability measure on
 $(\Omega, \F)$ defined by
 $$\begin{array}{l}
 \frac{d\PP^{u}}{d\PP} = \exp\{\int_0^T \sigma^{-1} (s,L.)
 f(s,L.,u_s) d\B_s - \frac{1}{2} \int_0^T |\sigma^{-1} (s,L.)
 f(s,L.,u_s)|^2 ds \}.\end{array}
 $$
 Thanks to Girsanov's Theorem (see e.g. Revuz and Yor (1991)),
 for every $u \in \V$  the process $\B^{u} := \left(\B_t - \int_0^t
 \sigma^{-1} (s,L_.) f(s,L_.,u_s) ds\right)_{0\le t \leq T}$ is a Brownian motion
 on ($\Omega,{\F},\PP^{u})$, and $L$ is a weak solution for the
 following functional differential equation.
 \begin{eqnarray}\label{app2}
 \left\{
 \begin{array}{ll}dL_t = f(t,L.,u_t) dt + \sigma(t,L.) d\B^{u}_t, \, \,\, \, 0< t \leq T, \\
 L_0 = x.
 \end{array}
 \right.
 \end{eqnarray}
 Under $\PP^u$, the process $L$ represents the evolution of the system when
 controlled by $(u_t)_{0\le t\leq T}$
 but not subject to impulses. Next, for a strategy $\delta=(\tau_n,\xi_n)_{n\geq 1}\in \A$, we denote by
$(L_t^\delta)_{0\le t \leq T}$ the process defined by
$$\begin{array}{ll}
L_t^\delta &= L_t + \sum_{n \geq 1} \xi_n
\ind_{[\tau_n < t]}\\
{}&=x+\int_0^tf(s,L.,u_s) ds + \int_0^t\sigma(s,L.) d\B^{u}_s+
\sum_{n \geq 1} \xi_n \ind_{[\tau_n < t]}.\end{array}
$$
Under $\PP^u$, the process $L^\delta$ stands for the evolution of
the system when controlled by $(u_t)_{0\le t\leq T}$ and subject to the
impulse strategy $\delta$. Note that the control and impulses are
interconnected. The reward function associated with the
 pair $(\delta,u)$ is
 \begin{equation}\label{reward-combi-u}
 J(\delta,u)=\E^u[\integ{0}{T}h(s,L^\delta, u_s)ds-\sum_{n\geq
 1}\psi(\xi_n) \ind_{[\tau_n<T]}],
 \end{equation}
  where, $\E^u$ is the expectation with respect to the probability measure
  $\PP^u$. With, $\xi_0=0$ and $\tau_0=0$, we have
$$
\int_{0}^{T}h(s,L^\delta, u_s)ds=\sum_{n\geq
0}\int_{\tau_n}^{\tau_{n+1}}h(s,L+\xi_1+\ldots+\xi_n, u_s)ds.
$$
The objective is to find a
  pair $(\delta^*,u^*)$ such that
$$
  J(\delta^*,u^*)=\sup_{(\delta,u)\in \A \times \cal V}
  J(\delta,u).
$$
Next let $H$ be the Hamiltonian associated with the control problem,
i.e., the function which with $(t,w, z,u) \in
 [0,T]\times {\cal C} \times {\mathbb{R}}^{d} \times \cal V$ associates
 $H(t,w,z,u) = z \sigma^{-1} (t,w)f(t,w,u)+h(t,w,u)$. The function
 $H$ is Lipschitz $w.r.t.$ $z$ uniformly in $(t,w,u)$ and through Bene\v s Selection Lemma (cf. Bene\v s (1970), Lemma 1), there
 exists a $\mathbb{G}\otimes {\cal{B}}(\R^d)-$measurable function
 with values in $\cal V$ such that for any $(t,w,z) \in
 [0,T] \times \R^{d+d}$,
 \begin{equation}\label{benes}
 H^*(t,w,z):=\sup_{u \in \cal V} H(t,w,z,u)=H(t,w,z,u^*(t,w,z)).
  \end{equation}
Moreover, the function $H^*$ is Lipschitz in $z$ uniformly w.r.t.
$(t,w)$ as a supremum over $u\in \cal V$ of functions uniformly
Lipschitz w.r.t. $(t,w,u)$.
\bigskip

For any stopping time $\nu \in \T$, and any ${\cal
F}_{\nu}$-measurable random variable $\xi$, let \\
$(Y^{n}(\xi,\nu),\;Z^{n}(\xi,\nu),\;K^{n}(\xi,\nu))_{n \geq 0}$ be
the sequence of processes defined as follows.
\begin{equation} \label{y0udef1}
Y^{0}_t(\nu,\xi)= \integ{t}{T}
H^*(s,L.(\omega)+\xi,Z_s^0(\xi,\nu))\ind_{[s\geq \nu]}ds
-\integ{t}{T} Z^0_s(\xi,\nu)d\B_s, \,\,\, 0\le t\leq T,
\end{equation}
and, for any  $n\geq 1$,
\begin{equation} \label{y0udef2}
\left\{
\begin{array}{l}
(Y^{n}(\nu,\xi),Z^{n}(\nu,\xi),K^{n}(\nu,\xi))\in {\cal S}_c^2\times {\H}^{2,d}\times \S_{c,i}^2 \\
Y^{n}_t(\nu,\xi)=\integ{t}{T}H^*(s,L.+\xi,Z_s^{n}(\xi,\nu))\ind_{[s\geq
\nu]}ds +K^{n}_T(\nu,\xi)- K^{n}_t(\nu,\xi)-\integ{t}{T}
Z^{n}_s(\nu,\xi)d\B_s,\,\,\, t\leq T,\\
Y^{n}_t(\nu,\xi)\geq O^{n}_t(\nu,\xi):=
\ds{\max_{\gb \in U}}\:(-\psi(\gb)+Y^{n-1}_t(\nu,\xi+\gb)), \,\,\, t \leq T,\\
\integ{0}{T}(Y^{n}_t(\nu,\xi)-O^{n}_t(\nu,\xi))dK^{n}_t(\nu,\xi)=0.
\end{array}
\right.
\end{equation}
We can easily see by induction that for any $n \geq 0$, the
processes $Y^{n}(\xi,\nu),Z^{n}(\xi,\nu)$ and $K^{n}(\xi,\nu)$ are
well defined, since $H^*$ is Lipschitz in $z$ and $U$ is finite. In
addition, the process $Y^{n}(\xi,\nu)$ is continuous, since
$\max_{\gb \in U}(-\psi(\gb)+Y^{n-1}_T(\nu,\xi+\gb))< 0$. Next,
in view of  Proposition \ref{thcompa}, it holds that, for
any $n\geq 0$, for any $\nu$ and $\xi$, $Y^{n}(\xi,\nu)\leq
Y^{n+1}(\xi,\nu)$ since $Y^{0}(\xi,\nu)\leq Y^{1}(\xi,\nu)$.

\medskip\no
Now, according to (\ref{y0udef1}) and (\ref{y0udef2}), there are
controls $u^n\in \V$ such that:
\begin{eqnarray} \label{snu0} Y_t^{0}(\nu,\xi) =
\E^ {u^0}[\int_t^T h(s,L.+\xi,u^0_s) \ind_{[s \geq \nu]} ds|\F_t],\;\,\,
  0\le t \leq T,
\end{eqnarray}
and,  for any $n \geq 1$,
\begin{eqnarray} \label{snnu}
Y_t^{n}(\nu,\xi) = \mbox{esssup}_{\tau \in \T_t}
\E^{u^n}[\int_t^\tau h(s,L.+\xi,u^n_s) \ind_{[s \geq \nu]}
ds+\ind_{[\tau <T]}O_\tau^{n} (\nu,\xi)|\F_t],\; t \leq T.
\end{eqnarray}
The last inequality is valid since $K^n(\nu,\xi)$ is non-decreasing
and $Y_t^{n}(\nu,\xi)\geq \ind_{[\tau <T]}O_\tau^{n} (\nu,\xi)$.
Therefore, $Y_t^{n}(\nu,\xi)$ is greater than the expression inside
the ess $\sup$. On the other hand,  there is equality when
$\tau=\inf\{s\ge t,\, K^n_s(\nu,\xi)-K^n_t(\nu,\xi)>0\}\wedge T$.

\no Now, by induction, as in the proof of Proposition \ref{propo1n}, we
obtain that, for any $n \geq 0$, $\tau$ a stopping time and any $\F_\tau$-measurable r.v. $\xi$, the process $Y^{n}(\nu,\xi)$ satisfies the following
property:
$$
0\leq Y_t^{n}(\nu,\xi) \leq \gamma (T-t),\;\:\,\, t \leq T,
$$
where, $\gamma$ is the constant of boundedness of $h$. Therefore,
using Proposition \ref{propopeng}, there
exists a \cadlag process $(Y_t^{*}(\nu,\xi))_{t \leq T}$ limit
of the increasing sequence $(Y^{n}(\nu,\xi))_{n\geq 0}$ as $n
\rightarrow \infty$. Moreover we have
$$
0\leq Y_t^*(\nu,\xi)\leq
\gamma (T-t), \quad t\leq T.
$$
In the next proposition, we give a characterization of $Y^*(\nu,\xi)$.

\begin{propo} \label{porop21}The process $Y^*(\nu,\xi)$ is
  continuous. Moreover, there exist processes
$Z^*(\nu,\xi)\in {\cal H}^{2,d}$
and $K^*(\nu,\xi)\in {\cal S}^2_{ci}$ such that, for all $t\leq T$,
\begin{equation} \label{eqconty2}\left\{
\begin{array}{lll}
Y^*_t(\nu,\xi)= \integ{t}{T}H^*(s,L.+\xi,Z_s^*(\xi,\nu))\ind_{[s\geq
\nu]}ds +K^{*}_T(\nu,\xi)- K^{*}_t(\nu,\xi)-\integ{t}{T}
Z^{*}_s(\nu,\xi)d\B_s,\\
Y^{*}_t(\nu,\xi)\geq O_t(\nu,\xi):= \ds{\max_{\gb \in
U}}\:(-\psi(\gb)+Y^{*}_t(\nu,\xi+\gb))\\
\integ{0}{T}(Y^{*}_t(\nu,\xi)-O_t(\nu,\xi))dK^{*}_t(\nu,\xi)=0.\end{array}\right.
\end{equation}
Furthermore, for any pair $(\nu,\xi)$ and any stopping time $\nu'\geq
\nu$, we have $Y^*_\nu(\nu,\xi)=Y_{\nu'}^*(\nu,\xi)$.
\end{propo}
\noindent\emph{Proof.} Thanks to Proposition
\ref{propopeng}, there exists a process $Z^*(\nu,\xi)\in {\cal
H}^{2,d}$ such that, for any $p\in [1,2)$, the sequence
$(Z^n(\nu,\xi))_{n\geq 0}$ converges to $Z^*(\nu,\xi)$ in ${\cal
H}^{p,d}$. This convergence holds also weakly in ${\cal H}^{2,d}$.
Additionally, there exists an increasing process $K^*(\nu,\xi)\in
{\cal S}^2_{i}$ such that for any stopping time $\tau$ the sequence
$(K^n_\tau(\nu,\xi))_{n\geq 0}$ converges to $K^*_\tau(\nu,\xi)$ in
$L^p(dP)$. Therefore, we have
\begin{equation}\label{eqqueq}\left\{\begin{array}{l}
Y^*_t(\nu,\xi)= \int_{t}^{T}H^*(s,L.+\xi,Z_s^*(\xi,\nu))\ind_{[s\geq
\nu]}ds +K^{*}_T(\nu,\xi)- K^{*}_t(\nu,\xi)-\int_{t}^{T}
Z^{*}_s(\nu,\xi)d\B_s,\\ Y^{*}_t(\nu,\xi)\geq
O_t(\nu,\xi):= \ds{\max_{\gb \in
U}}\:(-\psi(\gb)+Y^{*}_t(\nu,\xi+\gb)),\,\,\,\, 0\leq t\leq
T.\end{array}\right.
\end{equation}
The last inequality is valid, since $U$ is finite.

\no
Next, for $t\leq T$, let us set
\begin{equation}\label{eqr}R_t=\mbox{ess sup}_{\tau \in \T_t} \E[\int_t^\tau
H^*(s,L.+\xi,Z_s^*(\xi,\nu))\ind_{[s\geq \nu]}ds+\ind_{[\tau
<T]}O_\tau^* (\nu,\xi)|\F_t].\end{equation} Using Characterization
(\ref{snell}) of $(R_t)_{0\le t\leq T}$ as a solution of a BSDE yields
that, in using the comparison result (Proposition \ref{thcompa}),
for any $t\leq T$, $R_t\geq Y_t^{n}(\nu,\xi)$ and then $R_t\geq
Y_t^*(\nu,\xi)$. On the other hand, a result by Peng and Xu (2005)
implies that $(R_t)_{0\le t\leq T}$ is the smallest
$H^*(s,L.+\xi,z)\ind_{[s\geq \nu]}$-supermartingale which dominates
$O_t(\nu,\xi):= \ds{\max_{\gb \in
U}}\:(-\psi(\gb)+Y^{*}_t(\nu,\xi+\gb))$. But, by (\ref{eqqueq}), the
process $Y^{*}(\nu,\xi)$ is a $H^*(s,L.+\xi,z)\ind_{[s\geq
\nu]}$-supermartingale such that $Y^{*}_t(\nu,\xi)\geq
O_t(\nu,\xi):= \ds{\max_{\gb \in
U}}\:(-\psi(\gb)+Y^{*}_t(\nu,\xi+\gb))$. Thus, $Y^{*}_t(\nu,\xi)\geq
R_t$, for any $t\leq T$. Finally, since  both processes are \cadlag,
then P-$a.s.$, $R=Y^{*}(\nu,\xi)$. This means that $Y^{*}(\nu,\xi)$
is equal to the second term in (\ref{eqr}). Now, using the
characterization of Theorem \ref{thexist}, it holds that
$Y^*(\nu,\xi)$ and, $Z^*(\nu,\xi)$ and $K^*(\nu,\xi)$ satisfy
(\ref{eqconty2}). The continuity of $Y^*(\nu,\xi)$ is obtained in a
similar fashion as in Theorem \ref{veri-impulse} since $U$ is
finite.

\medskip\no Now, if $\nu'\geq \nu$ then thanks to uniqueness result we have, for
any $n\geq 0$, $Y^n_\nu(\nu,\xi)=Y_{\nu'}^n(\nu,\xi)$, and then it is
enough to take the limit as $n\to \infty$. $\Box$
\medskip

\no In the same way as previously, for any admissible control $u\in \V$,
a stopping time $\nu$, an  ${\cal F}_\nu$-measurable $r.v.$ $\xi$ and
$n\geq 0$, let us consider the sequence of processes defined
recursively by
\begin{equation} \label{y0udef}
Y^{u,0}_t(\nu,\xi)= \integ{t}{T}
H(s,L.(\omega)+\xi,Z_s^{u,0}(\xi,\nu),u_s)\ind_{[s\geq \nu]}ds
-\integ{t}{T} Z^{u,0}_s(\xi,\nu)d\B_s, \,\, t\leq T
\end{equation}
and, for any  $n\geq 1$,
\begin{equation} \label{y0udef}
\left\{
\begin{array}{l}
(Y^{u,n}(\nu,\xi),Z^{u,n}(\nu,\xi),K^{u,n}(\nu,\xi))\in {\cal S}_c^2\times {\H}^{2,d}\times \S_{c,i}^2 \\
Y^{u,n}_t(\nu,\xi)=\int_{t}^{T}H(s,L.+\xi,Z_s^{u,n}(\xi,\nu),u_s)\ind_{[s\geq
\nu]}ds +\\ \qquad \qquad\qquad \qquad\qquad\qquad
K^{u,n}_T(\nu,\xi)- K^{u,n}_t(\nu,\xi)-\int_{t}^{T}
Z^{u,n}_s(\nu,\xi)d\B_s,\\
Y^{u,n}_t(\nu,\xi)\geq O^{u,n}_t(\nu,\xi):=
\ds{\max_{\gb \in U}}\:(-\psi(\gb)+Y^{u,n-1}_t(\nu,\xi+\gb))\\
\mbox{ and
}\int_{0}^{T}(Y^{u,n}_t(\nu,\xi)-O^{u,n}_t(\nu,\xi))dK^{u,n}_t(\nu,\xi)=0.
\end{array}
\right.
\end{equation}
As above, the sequence of processes $(Y^{u,n}(\nu,\xi))_{n\geq
0}$ is increasing and converges to a \cadlag process $Y^u(\nu,\xi)$
which satisfies $0\leq Y_t^u(\nu,\xi)\leq \gamma (T-t)$, for any
$t\leq T$. We also have the following
\begin{propo} The process $Y^u(\nu,\xi)$ is continuous. Furthermore, there exist two
processes \\ $(Z^u(\nu,\xi),K^u(\nu,\xi))\in {\cal H}^{2,d}\times
{\cal S}^2_{ci}$ such that, for all $t\leq T$,
\begin{equation} \label{eqconty}\left\{
\begin{array}{l}
Y^u_t(\nu,\xi)=
\int_{t}^{T}H(s,L.+\xi,Z_s^u(\xi,\nu),u_s)\ind_{[s\geq \nu]}ds
+K^{u}_T(\nu,\xi)- K^{u}_t(\nu,\xi)-\int_{t}^{T}
Z^{u}_s(\nu,\xi)d\B_s,\\
Y^{u}_t(\nu,\xi)\geq O_t(\nu,\xi):= \ds{\max_{\gb \in
U}}\:(-\psi(\gb)+Y^{u}_t(\nu,\xi+\gb)), \\
\int_{0}^{T}(Y^{u}_t(\nu,\xi)-O_t(\nu,\xi))dK^{u}_t(\nu,\xi)=0.\end{array}\right.
\end{equation}
Moreover, we have
$$
Y^u_0(0,0)=\sup_{\delta \in {\cal
A}}J(u,\delta).
$$
\end{propo}

\noindent\emph{Proof.} The proof of the two first claims is the same
as the one of Proposition \ref{porop21}. It remains
to show the last one. Indeed, since the triple
$(Y^u(\nu,\xi), (Z^u(\nu,\xi),K^u(\nu,\xi))$ satisfies
\begin{equation} \label{eqconty}\left\{
\begin{array}{l}
Y^u_t(\nu,\xi)= \int_{t}^{T}h(s,L.+\xi,u_s)\ind_{[s\geq \nu]}ds
+K^{u}_T(\nu,\xi)- K^{u}_t(\nu,\xi)-\int_{t}^{T}
Z^{u}_s(\nu,\xi)d\B^u_s,\,\,\,\, t\leq T\\
Y^{u}_t(\nu,\xi)\geq O_t(\nu,\xi):= \ds{\max_{\gb \in
U}}\:(-\psi(\gb)+Y^{u}_t(\nu,\xi+\gb))\mbox{ and }\\
\int_{0}^{T}(Y^{u}_t(\nu,\xi)-O_t(\nu,\xi))dK^{u}_t(\nu,\xi)=0.\end{array}\right.
\end{equation}
it follows, as in Theorem \ref{veri-impulse}, that
$Y^u_0(0,0)=\sup_{\delta \in {\cal A}}J(u,\delta).$ $\Box$
\medskip

We give now the main result of this section.
\begin{thm}
There exist a control $u^* \in \V$ and a strategy $\delta^* =
(\tau_n^*,\gb_n^*)_{n \geq 0} \in \A$ such that
$$
J(\delta^*,u^*)=\sup_{(\delta,u)\in \A \times \cal V}
  J(\delta,u).
$$
In addition,
$$
Y_0^{u^*}(0,0) = J(\delta^*,u^*).
$$
\end{thm}
\no\emph{Proof: }Let $u \in \V$,  then through the definitions of
$Y^u(\nu,\xi)$ and $Y^*(\nu,\xi)$ it holds true that
$Y^*(\nu,\xi)\geq Y^u(\nu,\xi)$ since, in using the comparison
result of Proposition \ref{thcompa} and an induction argument, we have
$Y^{*,n}(\nu,\xi)\geq Y^{u,n}(\nu,\xi)$, for any $n\geq 0$.
Hence, we have
$$
Y_0^*(0,0)\geq Y_0^u(0,0) =\sup_{\delta \in {\cal A}}J(u,\delta),
$$
and then
$$
Y_0^*(0,0)\geq \sup_{(\delta,u)\in \A \times \cal V}
  J(\delta,u)\geq \sup_{u\in {\cal V}}\sup_{\delta \in {\cal
A}}J(u,\delta).
$$
Now, let $u^*$ and $\delta^*$ be defined as follows.
$$
\begin{array}{l}
\tau_1^{*} = inf\{s \geq 0;\; O_s (0,0)= Y_s^*(0,0)\}\wedge T,\\
-\psi(\beta_1^*) + Y^*_{\tau_{1}^*}
(\tau_{1}^*,\beta_1^*)=\max_{\beta \in U} \{-\psi(\beta) +
Y^*_{\tau_{1}^*} (0,\beta)\}=O_{\tau^*_1}(0,0),\\
u^*_t\1_{[t\leq \tau^*_1]}=u^*(t,L., Z^*_t(0,0))
\end{array}
$$
and, for $n\ge 2$,
$$
\begin{array}{llll}
\tau_n^{*} = inf\{s \geq \tau_{n-1}^{*}, \,Y_s^*
(\tau_{n-1}^{*},\beta_1^* +\ldots+\beta_{n-1}^*) =O_s
(\tau_{n-1}^{*},\beta_1^* +\ldots+\beta_{n-1}^*)\}\wedge
T,\\
-\psi(\beta_n^{*}) + Y_{\tau_{n}^{*}} (\tau_{n}^{*},\beta_1^{*}
+\ldots+\beta_{n-1}^{*}+\beta_n^{*})=\max_{\beta \in U} \{-\psi(\beta)
+ Y^*_{\tau_{n}^{*}} (\tau_{n-1}^{*},\beta_1^{*}
+\ldots+\beta_{n-1}^{*}+\beta)\}\\
\qquad\qquad\qquad\qquad\qquad\qquad\qquad\qquad\quad=O_{\tau_{n}^{*}}(\tau_{n-1}^{*},\beta_1^{*} +\ldots+\beta_{n-1}^{*})
\\
\mbox{and } u^*_t\1_{[\tau^*_{n-1},\tau^*_n]}(t)=u^*(t,L.+\beta_1^*
+\ldots+\beta_{n-1}^*, Z^*_t(\tau_{n-1}^{*},\beta_1^*
+\ldots+\beta_{n-1}^*)).
\end{array}
$$
Therefore,
\begin{eqnarray*}
Y^{*}_0(0,0)&=& \E^{u^*}[\int_0^{\tau^{*}_1} h(s,L,u^*_s)ds +
O_{\tau^*_1}(0,0)\1_{[\tau^*_1<T]}]
\end{eqnarray*}
and as $O_{\tau^*_1}(0,0)=-\psi(\beta_1^*) + Y^*_{\tau_{1}^*}
(\tau_{1}^*,\beta_1^*)$ then
\begin{eqnarray*}
Y^{*}_0(0,0)&=& \E^{u^*}[\int_0^{\tau^{*}_1} h(s,L,u^*_s)ds +
(-\psi(\beta_1^*) + Y^*_{\tau_{1}^*} (\tau^*_1,
\beta_1^*))\1_{[\tau^*_1<T]}].
\end{eqnarray*}
But,
$$\begin{array}{ll}
Y^*_{\tau_{1}^*} (\tau^*_1,\beta_1^*)&=Y^*_{\tau_{2}^*}
(\tau^*_1,\beta_1^*)+
\integ{\tau_{1}^*}{\tau_{2}^*}h(s,L.+\beta_1^*,u^*_s)ds
-\integ{\tau_{1}^*}{\tau_{2}^*}
Z^{*}_s(\tau^*_1,\beta_1^*)d\B^{u^*}_s\\{}&= \E^{u^*}[
Y^*_{\tau_{2}^*} (\tau^*_1,\beta_1^*)+
\integ{\tau_{1}^*}{\tau_{2}^*}h(s,L.+\beta_1^*,u^*_s)ds|{\cal
F}_{\tau_1^*}] .\end{array}
$$
Plugging the last quantity in the previous equality to obtain
\begin{eqnarray*}
Y^{*}_0(0,0)&=& \E^{u^*}[\int_0^{\tau^{*}_1}
h(s,L,u^*_s)ds+\int_{\tau^{*}_1}^{\tau^{*}_2}
h(s,L+\beta_1^*,u^*_s)ds -\psi(\beta_1^*)\1_{[\tau^*_1<T]} +
Y^*_{\tau_{2}^*} (\tau^*_1, \beta_1^*)\1_{[\tau^*_2<T]}]\\
{}&=& \E^{u^*}[\int_0^{\tau^{*}_2} h(s,L^{\delta^*},u^*_s)ds
-\psi(\beta_1^*)\1_{[\tau^*_1<T]} + Y^*_{\tau_{2}^*} (\tau^*_1,
\beta_1^*)]\\{}&=&\E^{u^*}[\int_0^{\tau^{*}_2}
h(s,L^{\delta^*},u^*_s)ds -\psi(\beta_1^*)\1_{[\tau^*_1<T]} +
O_{\tau_{2}^*} (\tau^*_1, \beta_1^*)\1_{[\tau^*_2<T]} ],
\end{eqnarray*}
since $Y^*_{\tau_{2}^*} (\tau^*_1,\beta_1^*)=Y^*_{\tau_{2}^*}
(\tau^*_1,\beta_1^*)\1_{[\tau^*_2<T]}$,  $[\tau^*_2<T]\subset
[\tau^*_1<T]$ and finally \\ $Y^*_{\tau_{2}^*} (\tau^*_1,
\beta_1^*)=O^*_{\tau_{2}^*} (\tau^*_1, \beta_1^*)\1_{[\tau^*_2<T]}$.

\no Repeating now this reasoning as many times as necessary to obtain,
for all $n\geq 1$,
\begin{eqnarray*}
Y^{*}_0(0,0)=\E^{u^*}[\int_0^{\tau^{*}_n} h(s,L^{\delta^*},u^*_s)ds
-\sum_{k=1,n}\psi(\beta_k^*)\1_{[\tau^*_k<T]} + O_{\tau_{n+1}^*}
(\tau^*_n,\beta_1^*+...+ \beta_n^*)\1_{[\tau^*_{n+1}<T]} ].
\end{eqnarray*}
This property implies first that the strategy $\delta^*$ is finite
since $Y^*(0,0)$ is a real constant. On the other hand taking the
limit as $n\to \infty$ to obtain:
$$
Y^*(0,0)=J(u^*,\delta^*).
$$
Thus,
$$
Y^*(0,0)=J(\delta^*,u^*) = \sup_{u \in \V} \;\sup_{\delta \in \A}\;
J(\delta,u),
$$
and the proof is complete. $\Box$
 \section*{Appendix}
 Let $\theta$ (resp. $\pi$) be  the optional (resp. predictable)
 tribe on $\Omega \times [0,T]$, i.e., the tribe  generated by the
 \cadlag and $\F_t$-adapted processes $X=(X_t)_{0\le t\leq T}$ (resp.
  the left continuous and $\F_t$-adapted  processes
 $Y=(Y_t)_{t\leq T}$ ).
 \medskip

 \begin{definition}\label{defa1}
  A measurable process $U=(U)_{t\leq T}$ is said to be of class
 [D] if the set of random variables $\{U_\tau, \tau \in {\cal T}\}$
 is uniformly integrable.
 \end{definition}
 \begin{propo}
 Let $U=(U)_{t\leq T}$ be an optional process which is of class [D]
 and $N= (N_t)_{t \leq T}$ the Snell envelope of $U$ defined by:
 $$
 N_t = \mbox{ess sup}_{\tau \in \T_t} \E[U_\tau |\F_t],\; t \leq T.
 $$
  If $U$ is right upper
 semi-continuous, then the process $N$ is continuous.
 \end{propo}
 \begin{propo}\label{prapp2}
  Let $(U)_{t \leq T}$ be an optional process of class [D]
 and $N$ its Snell envelope. Then
 \begin{description}
 \item[($i$)]   there exist a martingale $M$ and two increasing,
 integrable and  right continuous processes $A$ and $B$ such that,
 \begin{equation}
  N_t = M_t - A_{t} - B_t,\quad 0\le t\le T.
 \end{equation}
 The process $A$ is optional and continuous, and $B$ is
 predictable, i.e., $\pi-$measurable and purely discontinuous. This
 decomposition is unique. In addition for any $t \leq T$ we have:
 \begin{equation}\label{eqann}
 \{\Delta_t B >0\} \subset \{U_{t-} = N_{t-}\}
 \end{equation}
 and
 \begin{equation}
 \Delta_t B  =(U_{t-} -N_{t-})^+ \ind_{[\Delta_t U <0]}.
 \end{equation}
  \item[($ii$)] If  $Y \in \S^2$ and $M$ is a continuous
 martingale with respect to $\F$, then the processes $A$ and $B$
 are also in $\S^2$.
 \end{description}
 \end{propo}

 {\small{
 }}
 \end{document}